# A Comparative Analysis of the Relative Efficacy of Vector-Control Strategies against Dengue Fever

*Marcos Amaku[1], Francisco Antonio Bezerra Coutinho[2], Silvia Martorano Raimundo[2], Luis Fernandez Lopez[2,3], Marcelo Nascimento Burattini[2], Eduardo Massad[2,4,*]*

[1]School of Veterinary Medicine, University of São Paulo, Av.Prof. Orlando Marques de Paiva, 87 – Cidade Universitária, São Paulo/SP – CEP 05508 270, Brazil

[2]School of Medicine, University of Sao Paulo and LIM 01-HCFMUSP, Av. Dr. Arnaldo 455, São Paulo/SP – CEP 01246-903, Brazil

[3] CIARA - Florida International University, Miami, USA

[4]London School of Hygiene and Tropical Medicine, London University, Keppel Street, London W1C 7HT, UK

[*]Corresponding author: edmassad@usp.br

**Short title: Relative Efficacy of Dengue Vectors Control**

**Abstract**

**Background.** Dengue is considered one of the most important vector-borne infection, affecting almost half the world population with 50 to 100 millions cases every year. In this paper, we present one of the simplest models that can encapsulate all the important variables related to vector control in dengue fever.

**Methodology.** The model considers the human population, the adult mosquito population and the population of immature stages, which includes eggs, larvae and pupae. The model also considers the vertical transmission of dengue in the mosquitoes and the seasonal variation in the mosquito population. From this basic model describing the dynamics of dengue infection, we deduce thresholds for avoiding the introduction of the disease and for the elimination of the disease. In particular, we deduce a Basic Reproduction Number for dengue that includes parameters related to the immature stages of the mosquito. By neglecting seasonal variation, we calculate the equilibrium values of the model's variables. We also present a sensitivity analysis of the impact of four vector-control strategies on the Basic Reproduction Number, on the Force of Infection and on the human prevalence of dengue. Each of the strategies was studied separately from the others.

**Principal Findings.** The analysis presented allows us to conclude that of the available vector control strategies, adulticide application is the most effective, followed by the reduction of the exposure to mosquito bites, locating and destroying breeding places and, finally, larvicides.

**Significance.** Current vector-control methods are concentrated in mechanical destruction of mosquitoes' breeding places. Our results suggest that reducing the contact between vector and hosts (biting rates) are as efficient as the logistically difficult but very efficient adult mosquito's control.

**Author summary**

Dengue is a viral disease that affects almost half of the world population. There is no specific treatment and the vaccine is still in its first trials and will not be available for the next three or four years. Controlling the disease is therefore restricted to reduce the number of its vector, mosquitoes from the genus *Aedes.* The available vector-control strategies are the mechanical destruction of the mosquito breeding places (the vector breeds in urban or peri-urban environment), larvicides and adulticides.

Here we propose a mathematical model that captures the essence of dengue transmission, from which we derive the main parameters related to the intensity of dengue transmission, namely the Basic Reproduction Number, on the Force of Infection, and on the human prevalence of dengue. We analyze which the parameters of the model these quantities are most sensitive to. We also analyzed the model's sensitivity to the mosquitoes' biting rate and showed that reducing this parameter with repellents and mosquitoes' shields (clothes impregnated with insecticides), along with the increase in the adult mosquito mortality rate by the use of insecticides are the most effective control strategies against dengue.

## Introduction

The global expansion of dengue fever is a matter of great concern to public health authorities around the world [1]. In terms of the population at risk, dengue is considered the most important vector-borne disease worldwide. [2,3]. It is estimated that approximately 3.6 billion people, one-half of the world's population, live in parts of the world affected by dengue [4-6], and 120 million people are expected to travel to dengue-affected areas every year [7]. Between 50 and 100 million people are infected each year [8], and the World Health Organization states that the number is rising due to human population growth and the increased spread of vector mosquitoes due to climate change [9]. Recent studies suggest that the figures are much higher [10], with as many as 230 million infections, tens of millions of cases of dengue fever (DF) and millions of cases of dengue hemorrhagic fever DHF [6,11,12]. The number of disability-adjusted life years (DALYs) worldwide is estimated to range between 528 and 621 per million population [10, 13], and the total cost of dengue cases in the affected areas of the world may be approximately 2 billion dollars annually [8].

Dengue viruses are transmitted by mosquitoes of the genus *Aedes*, subgenus *Stegomyia* [10]. The principal vector, *Aedes Stegomyia aegypti*, is now well established in much of the tropical and subtropical world, particularly in urban areas. It is a domestic species, highly susceptible to dengue virus infection, feeding preferentially on human blood during the daytime and often taking multiple blood meals during a single gonotrophic cycle [13]. It typically breeds in clean stagnant water in artificial containers and is, therefore, well adapted to urban life. A second species, *Aedes Stegoymyia albopictus*, is generally considered less effective as an epidemic vector because, unlike *A. aegypti*, it feeds on many animals other than humans and is less strongly associated with the domestic environment [14].

Several reasons have been proposed for the dramatic global emergence of dengue as a major public health problem. Major global demographic changes have occurred, the most important of which have been uncontrolled urbanization and concurrent population growth. The public health infrastructure of many of the affected countries has deteriorated. Increases in international travel provide an efficient mechanism for the human transport of dengue viruses between urban centers, resulting in the frequent exchange of dengue viruses. Climatic changes influence the mosquito's survival and proliferation [15]. Finally, effective mosquito control is virtually nonexistent in many dengue-endemic countries [16, 17].

99    Essentially, the control of dengue has been based on three strategies [18]: source
100   reduction (locating and destroying mosquitoes' breeding places), larvicides and ultra-
101   low volume (ULV) application of aerosol adulticides. The first two strategies have been
102   applied with varying degrees of success. However, there is still considerable
103   controversy over the efficacy of the current methods for controlling adult mosquitoes
104   [18]. At the time of the advent of DDT, *Aedes aegypti* was highly susceptible to this
105   agent [18]. The successful application of DDT resulted in the eradication of *Aedes*
106   *aegypti* from 22 countries in the Americas in 1962 and from all countries in the
107   Mediterranean region in 1972. However, the fate of DDT is well known. DDT was
108   abandoned due to the evolution of resistant insects and due to the environmental
109   impacts of the insecticide. Therefore, the control of dengue shifted to other approaches:
110   source reduction, larvicides and adulticides from other chemical families.

111   From a theoretical perspective, significant advances were made by Macdonald
112   [19], who proposed that the most effective control strategy against vector-borne
113   infections is to kill adult mosquitoes.

114   Recently, in a study for describing the dynamics of dengue, we showed that the
115   models describing infections transmitted by blood-sucking insects are indeed very
116   sensitive to the mosquitoes' mortality rate [20].

117   The current paper differs from the ones previously published [15, 20, 21] in the
118   following sense: In [21] the model's basic structure is presented, in particular it presents
119   a new seasonality factor. Thus, paper [21] was designed to test one hypothesis to
120   explain dengue's overwintering; in [20] the model presented in [21] was numerically
121   simulated in order to mimic Singapore data. In addition, an incomplete sensitivity
122   analysis was presented, which intended to show that killing adult mosquitoes was the
123   most effective strategy, as demonstrated numerically in that paper. The role of larvicide
124   as an important tool to avoid the resurgence of outbreaks was proposed based only in
125   numerical simulations. The paper by Massad et al. [15] is a review of the previous
126   papers and does not add anything new on control.

127   The current paper is an analysis of the basic model proposed in [21] and
128   numerically studied in [20].

129   In this paper, we present what we consider to be the simplest model that
130   encapsulates all the important variables related to dengue control, and we analyze four
131   control strategies for use against the vectors of dengue. All the relevant stages are
132   included and the ones not included (like larvae and pupae) can be trivially added to the

133 model and a complete analysis of the sensitivity of transmission to the parameters is
134 presented.

135

136

137 **Methods**

138 **The basic model**

139     The basic model that is used to calculate the efficiency of control strategies can
140 be found in [15, 20, 21].

141     The populations involved in the transmission are human hosts, mosquitoes and
142 their eggs. For the purposes of this paper, the term "eggs" also includes the intermediate
143 stages, such as larvae and pupae. Therefore, the population densities are divided into the
144 following compartments: susceptible humans denoted $S_H$; infected humans, $I_H$;
145 recovered (and immune) humans, $R_H$; total humans, $N_H$; susceptible mosquitoes, $S_M$;
146 infected and latent mosquitoes, $L_M$; infected and infectious mosquitoes, $I_M$; non-infected
147 eggs, $S_E$; and infected eggs, $I_E$. The variables appearing in the model are summarized in
148 Table 1.

149     The model is defined by the following equations:

150

$$\frac{dS_H}{dt} = -abI_M \frac{S_H}{N_H} - \mu_H S_H + r_H N_H \left(1 - \frac{N_H}{\kappa_H}\right)$$

$$\frac{dI_H}{dt} = abI_M \frac{S_H}{N_H} - (\mu_H + \alpha_H + \gamma_H)I_H$$

$$\frac{dR_H}{dt} = \gamma_H I_H - \mu_H R_H$$

$$\frac{dS_M}{dt} = pc_S(t)S_E - \mu_M S_M - acS_M \frac{I_H}{N_H}$$

$$\frac{dL_M}{dt} = acS_M \frac{I_H}{N_H} - \gamma_M L_M - \mu_M L_M$$

$$\frac{dI_M}{dt} = \gamma_M L_M - \mu_M I_M + pc_S(t)I_E$$

$$\frac{dS_E}{dt} = \left[r_M S_M + (1-g)r_M(I_M + L_M)\right]\left(1 - \frac{(S_E + I_E)}{\kappa_E}\right) - \mu_E S_E - pc_S(t)S_E$$

$$\frac{dI_E}{dt} = \left[gr_M(I_M + L_M)\right]\left(1 - \frac{(S_E + I_E)}{\kappa_E}\right) - \mu_E I_E - pc_S(t)I_E$$

$$N_H = S_H + I_H + R_H$$
$$N_M = S_M + L_M + I_M$$

151

$$N_E = S_E + L_E$$

(1)

152

153

154 where $c_s(t) = \left(d_1 - d_2 sin\left(2\pi ft + \phi\right)\right)$ is a factor mimicking seasonal influences in the

155 mosquito population [21,22].

156

157 *Remark: This model differs from the classical Ross-Macdonald model because the*

158 *extrinsic incubation period in the classical Ross-Macdonald model is assumed to last*

159 $\tau$ *days, whereas in model (1) we assumed an exponential distribution for the latency in*

160 *the mosquitoes. The classical Ross- Macdonald model can be obtained from system (1)*

161 *by replacing the fifth and sixth equations by*

162

163
$$\frac{dL_M}{dt} = acS_M \frac{I_H}{N_H} - \mu_M L_M - acS_M\left(t-\tau\right)\frac{I_H\left(t-\tau\right)}{N_H\left(t-\tau\right)}e^{-\mu_M \tau}$$

$$\frac{dI_M}{dt} = acS_M\left(t-\tau\right)\frac{I_H\left(t-\tau\right)}{N_H\left(t-\tau\right)}e^{-\mu_M \tau} - \mu_M I_M + pc_s(t)I_E$$

164

165 *where $\tau$ is the extrinsic incubation period and $\mu_M$ is the mosquito mortality rate. The*

166 *expressions developed below in this paper with equations (1) can be replaced by the*

167 *corresponding expressions of the classical Ross-Macdonald model described above by*

168 *replacing $\frac{\gamma_M}{\gamma_M + \mu_M}$ by $e^{-\mu_M \tau}$. $\gamma_M$ is related to $\tau$ by $\tau = \frac{1}{\mu_M}\ln\left[\frac{\gamma_M}{\gamma_M + \mu_M}\right]$.*

169

170 **Equilibrium densities in the absence of seasonality**

171 The equilibrium densities of model (1) can be calculated exactly in the case

172 where seasonality can be neglected, i.e., with $c_s(t) = c_s =$ constant.

173 We begin by examining the steady-state values with $\alpha_H = 0$, i.e., with no

174 disease-induced mortality in the human population. Because we set $\alpha_H = 0$, we denote

175 the model variables with a superscript zero. By setting the derivatives in system (1) and

176 $\alpha_H$ equal to zero, it is straightforward to solve the resulting system of nonlinear

177 equations. The results are:

178

179 $N_H^0 = \kappa_H\left(\frac{r_H - \mu_H}{r_H}\right)$ (2)

180

181 $$N_M = N_M^0 = \frac{pc_S}{\mu_M} \kappa_E \left[ 1 - \frac{(\mu_M)(\mu_E + pc_S)}{r_M \, pc_S} \right] \tag{3}$$

182

183 $$N_E = N_E^0 = \kappa_E \left[ 1 - \frac{(\mu_M)(\mu_E + pc_S)}{r_M \, pc_S} \right] \tag{4}$$

184

185 Note that $N_M$ and $N_E$ do not depend on the disease mortality in the human population,

186 i.e., they do not depend on $\alpha_H$.

187

188 $$I_H^0 = \frac{(\gamma_M + g\mu_M)a^2 bc N_M - N_H^0 (\mu_H + \gamma_H)(\mu_M + \gamma_M)\mu_M (1-g)}{(\gamma_M + g\mu_M)a^2 bc \frac{N_M}{N_H^0}\left(1 + \frac{\gamma_H}{\mu_H}\right) + ac(\mu_H + \gamma_H)(\mu_M + \gamma_M)} \tag{5}$$

189

190 $$R_H^0 = \frac{\gamma_H}{\mu_H} I_H^0 \tag{6}$$

191 $$S_H^0 = N_H^0 - I_H^0 - R_H^0 \tag{7}$$

192

193

194 $$S_M^0 = \frac{(1-g)r_M N_M (\kappa_E - N_E) pc_S}{\kappa_E \left( \mu_M + ac \frac{I_H^0}{N_H^0} \right)(\mu_E + pc_S) - gr_M pc_S (\kappa_E - N_E)} \tag{8}$$

195

196 $$I_M^0 = \frac{(\mu_H + \gamma_H) I_H^0}{ab\left(1 - \left(1 + \frac{\gamma_H}{\mu_H}\right)\frac{I_H^0}{N_H^0}\right)} \tag{9}$$

197

198 $$L_M^0 = \frac{ac \frac{I_H^0}{N_H^0} S_M^0}{\gamma_M + \mu_M} \tag{10}$$

199
200

201 $$S_E^0 = \frac{\left[ r_M S_M^0 + (1-g) r_M (N_M - S_M^0) \right](\kappa_E - N_E)}{\kappa_E (\mu_E + pc_S)} \tag{11}$$

202

203 $$I_E^0 = N_E^0 - S_E^0 \tag{12}$$

204
205
206

207  If $\alpha_H \neq 0$, the total numbers of mosquitoes and eggs do not change. The expression for

208  $N_H$ is complicated, but it is straightforward to calculate $I_H$ as a function of $N_H$ as

209  follows:

210

211
$$\frac{I_H}{N_H} = \frac{(\gamma_M + g\mu_M)a^2bc\dfrac{N_M}{N_H} - (\mu_H + \gamma_H + \alpha_H)(\mu_M + \gamma_M)\mu_M(1-g)}{(\gamma_M + g\mu_M)a^2bc\dfrac{N_M}{N_H}\left(1 + \dfrac{\gamma_H}{\mu_H}\right) + ac(\mu_H + \gamma_H + \alpha_H)(\mu_M + \gamma_M)} \tag{13}$$

212

213  Alternatively, we can write

214

215
$$\frac{I_H}{N_H} = -\frac{\mu_H}{\alpha_H} + \frac{r_H}{\alpha_H}\left(1 - \frac{N_H}{\kappa_H}\right) \tag{14}$$

216

217

218  If the disease induces mortality in the human population ($\alpha_H \neq 0$), $N_H$ depends on $\alpha_H$

219  and is specified by a somewhat complicated expression. We will first obtain an

220  expression for $N_H$ as a function of $\alpha_H$. This expression is based on perturbation theory.

221  The exact expression for $N_H$ is presented subsequently.

222

223  **Estimating $N_H$ by perturbation theory**

224  An expression for $N_H$ can be obtained with perturbation theory. First, we sum the

225  first three equations of system (1) to obtain

226

227
$$\frac{dN_H}{dt} = r_H N_H\left(1 - \frac{N_H}{\kappa_H}\right) - \mu_H N_H - \alpha_H I_H \tag{15,}$$

228

229  At equilibrium, this expression yields

230

231
$$r_H N_H\left(1 - \frac{N_H}{\kappa_H}\right) - \mu_H N_H - \alpha_H I_H = 0 \tag{16}$$

232

233  Next, we expand $N_H$ and $I_H$ in powers of $\alpha_H$:

234

235
$$N_H = N_H^0 + \alpha_H N_H^1 + \alpha_H^2 N_H^1 + \mathrm{O}(\alpha_H^3) \tag{17}$$

236

237
$$I_H = I_H^0 + \alpha_H I_H^1 + \alpha_H^2 I_H^1 + \mathrm{O}(\alpha_H^3) \tag{18}$$

238

239      Neglecting the higher-order terms (because $\alpha_H$ is assumed to be very small) in (17) and

240      (18) and substituting in (16), we obtain, after some algebraic manipulations:

241

242
$$N_H = N_H^0 - \frac{\alpha_H I_H^0}{r_H - \mu_H} \tag{19}$$

243

244      where $N_H^0$ and $I_H^0$ are given by equations (2) and (5).

245

246

247 **The exact calculation of** $N_H$

248      The value of $\alpha_H$ for dengue is such that an individual who is sick for five days

249      has a probability of dying of the order of 0.2%, i.e., a negligible impact on human

250      demography. However, although it is reasonable to neglect $\alpha_H$ for dengue, it is not

251      reasonable to do so for other vector-borne infections, such as yellow fever or malaria.

252      We therefore need the exact expression for $N_H$ given below.

253      First, we define:

254

255
$$\Gamma = (\gamma_M + g\mu_M)a^2 bc N_M \tag{20}$$

256

257      where $N_M$ is given by equation (3), and

258

259
$$\theta = (\mu_H + \gamma_H + \alpha_H)(\mu_M + \gamma_M) \tag{21}$$

260

261      Next, we define:

262

263
$$\Pi = acr_H\theta \tag{22},$$

264
$$\Theta = -\left[ ac\theta\kappa_H(r_H - \mu_H) - \Gamma r_H\left(1 + \frac{\gamma_H}{\mu_H}\right) + \theta\mu_M\alpha_H\kappa_H(1 - g) \right] \tag{23}$$

265      and

266
$$\Omega = -\Gamma\kappa_H(r_H - \mu_H)\left(1 + \frac{\gamma_H}{\mu_H}\right) + \Gamma\alpha_H\kappa_H \tag{24}$$

267

268      Finally,

269

270 $$N_H = \frac{-\Theta + \sqrt{\Theta^2 - 4\Pi\Omega}}{2\Pi}$$ (25)

271

272 This expression reduces to equation (2) if $\alpha_H = 0$.

273

274

275 **Sensitivity of the variables to the parameters**

276      If seasonality is neglected (i.e., $c_s(t) =$ constant), the variables attain steady

277 states, as we have shown above. To estimate the sensitivity of a model variable in

278 steady state, $V_i$, to a parameter $\theta_j$, we consider the relative variation in the parameter,

279 $\frac{\Delta\theta_j}{\theta_j}$. This variation will correspond to a variation $\frac{\Delta V_i}{V_i}$ in the model variable $V_i$ given

280 by:

281

282

283

284 $$\frac{\Delta V_i}{V_i} = \frac{\theta_j}{V_i} \frac{\left[V_i(\theta_j + \Delta\theta_j) - V_i(\theta_j)\right]}{\Delta\theta_j} \frac{\Delta\theta_j}{\theta_j}$$ (26).

285

286

287 This expression can be approximated by [23,24]:

288

289

290 $$\frac{\Delta V_i}{V_i} = \frac{\theta_j}{V_i} \frac{\partial V_i}{\partial\theta_j} \frac{\Delta\theta_j}{\theta_j} + \frac{1}{2!} \frac{\theta_j^2}{V_i^2} \frac{\partial^2 V_i}{\partial\theta_j^2} \left(\frac{\Delta\theta_j}{\theta_j}\right)^2 + \dots$$ (27)

291

292 Usually, the second- and higher-order terms can be neglected provided that the relative

293 variation in the parameter, $\frac{\Delta\theta_j}{\theta_j}$, is sufficiently small.

294

295 **The sensitivity of the Basic Reproduction Number to the model's parameters**

296      Linearizing the second, the fifth, the sixth and the eight equations of model (1)

297 around the trivial solution (no-infection), we obtain the threshold normally denoted $R_0$

298 in the literature (details can be found in [21]).

299

300 $$R_0 = \frac{a^2 bc\left(\overline{N}_M / \overline{N}_H\right)\left(g\mu_M + \gamma_M\right)}{(\mu_H + \alpha_H + \gamma_H)(\mu_M + \gamma_M)\mu_M(1-g)}$$ (28)

301

302 where $\overline{N}_M$ and $\overline{N}_H$ denote the density of mosquitoes and of humans in the absence of

303 disease, respectively. Note that if $g = 0$, i.e., no vertical transmission, the expression

304 (28) for $R_0$ reduces to the classical Macdonald equation [25]. As mentioned above, to

305 obtain the classical Macdonald equation we replace $\dfrac{\gamma_M}{\gamma_M + \mu_M}$ by $e^{-\mu_M \tau}$. The case of

306 $g \to 1$ will be examined in the Discussion section.

307       Alternatively, we can deduce a threshold, $T_h$, for the existence of endemic

308 equilibrium values for the human prevalence of the disease. This threshold is given by

309 equation (13):

310
311

312 $$\frac{I_H}{N_H} = \frac{(\gamma_M + g\mu_M)a^2 bc\dfrac{N_M}{N_H} - (\mu_H + \gamma_H + \alpha_H)(\mu_M + \gamma_M)\mu_M(1-g)}{(\gamma_M + g\mu_M)a^2 bc\dfrac{N_M}{N_H}\left(1 - \dfrac{\gamma_M}{\mu_H}\right) + ac(\mu_H + \gamma_H + \alpha_H)(\mu_M + \gamma_M)}$$

313
314
315 If
316
317

318 $$\frac{I_H}{N_H} \geq 0 \; ,$$

319
320
321 then an endemic state exists. For this outcome, it suffices that

322
323

324 $$(\gamma_M + g\mu_M)a^2 bc\frac{N_M}{N_H} - (\mu_H + \gamma_H + \alpha_H)(\mu_M + \gamma_M)\mu_M(1-g) \geq 0$$ (29)

325
326
327 or
328

329 $$T_h = \frac{a^2 bc\left(N_M / N_H\right)\left(g\mu_M + \gamma_M\right)}{(\mu_H + \alpha_H + \gamma_H)(\mu_M + \gamma_M)\mu_M(1-g)} \geq 1$$ ,

330
331

332     which coincides with expression (28) if $T_h \leq 1$ because then $N_M = \overline{N}_M$ and $N_H = \overline{N}_H$.

333     This result also holds if $\alpha_H = 0$, i.e., if the disease has no influence on the population

334     size. Note that in our model, because the disease has no influence on the size of the

335     mosquito population, $N_M = \overline{N}_M$ always holds.

336        We begin the sensitivity analysis by considering the impact of a form of control

337     of dengue vectors that is still unusual, namely, reducing the contact of the population

338     with mosquito bites. This form of control is represented by mosquito shields (repellent-

339     impregnated cloths), repellents and the use of bed-nets. The use of bed-nets is very

340     effective against malaria [26] because it reduces the amount of contact between the

341     anopheline vectors and susceptible humans, the biting rate parameter $a$ of model (1).

342     We are aware that this strategy is effective against Anopheles mosquitoes because these

343     vectors bite at twilight and early at night. In contrast, Aedes mosquitoes bite primarily

344     during the day. We include this analysis here for the sake of generality and also because

345     the use of repellents and mosquito shields can produce the same reduction in the biting

346     rate $a$ and can be applied against Aedes mosquitoes. The partial derivative of $R_0$ with

347     respect to $a$ is given by

348

349
$$\frac{\partial R_0}{\partial a} = \frac{R_0}{a} \left[ 2 - \frac{a}{N_H} \frac{\partial N_H}{\partial a} \right] \qquad (30)$$

350

351     Next, we analyze the impact of reducing the carrying capacity of the immature forms,

352     $\kappa_E$, on the magnitude of $R_0$. This reduction represents a component of the strategy of

353     mechanical control, i.e., the identification and destruction of the places where Aedes

354     mosquitoes breed. The partial derivative of $R_0$ with respect to $\kappa_E$ is given by

355

356
$$\frac{\partial R_0}{\partial \kappa_E} = R_0 \left\{ \frac{p c_S}{N_M \mu_M} \left[ 1 - \frac{\mu_M (\mu_E + p c_S)}{r_M \, p c_S} \right] - \frac{1}{N_H} \frac{\partial N_H}{\partial \kappa_E} \right\} \qquad (31)$$

357

358

359     The use of larvicides is assumed to increase the mortality rate of the larvae, $\mu_E$.

360     Therefore, the impact of such a strategy is a function of the partial derivative of $R_0$ with

361     respect to $\mu_E$, which is

362

$$\frac{\partial R_0}{\partial \mu_E} = -R_0 \left( \frac{\kappa_E}{r_M N_M} + \frac{1}{N_H} \frac{\partial N_H}{\partial \mu_E} \right) \tag{32}$$

364

Finally, we take the partial derivative of $R_0$ with respect to the mosquito mortality rate

$\mu_M$ to estimate the impact of the application of adulticides as a control strategy against

the dengue vectors. The result is given by

368

$$\frac{\partial R_0}{\partial \mu_M} = R_0 \left[ \frac{1}{\mu_M + \gamma_M} + \frac{1}{\mu_M (1-g)} - \frac{p c_S \kappa_E}{\mu_M^2 N_M} - \frac{1}{N_H} \frac{\partial N_H}{\partial \mu_M} \right] \tag{33}$$

370

Given these partial derivatives, we can calculate the sensitivity of $R_0$ to the four

parameters above and thereby estimate the relative efficiencies of the control strategies

for avoiding the introduction of dengue into a non-infected area. To perform these

calculations, we consider equation (27) for each of the parameters. For dengue, the last

term in equations (30)-(33), involving the derivative of $N_H$, is always very small

relative to the previous terms. The results of the sensitivity analysis, with parameters'

values as in Table 2, are shown in Table 3.

378

379

**The sensitivity of the Force of Infection and the human prevalence to the model's**

**parameters**

382

The concept of 'force of infection' for vector-borne infection first appears in the

seminal works of Ronald Ross [27], who termed it the effective inoculation rate and

denoted it as $h$, for 'dependent happening'. The concept was further elaborated by

George MacDonald [19] who, in a now-famous appendix to his paper 'The Analysis of

Equilibrium in Malaria', defined the inoculation rate as

388

$h = mabs \tag{34}$

390

391    where $m$ is the mosquito density relative to the human population ($\frac{N_M}{N_H}$ in our notation),

392    $a$ is the mosquito's daily rate of biting, $b$ is the probability of infection from

393    mosquitoes to humans and $s$ is a quantity that Macdonald termed the 'Sporozoite Rate',

394    i.e., the prevalence of infection in the mosquitoes ($\frac{I_M}{N_M}$ in our notation). Note that

395    equation (34) is now expressed as

396    $$\lambda = ab\frac{I_M}{N_H} \tag{35}$$

397

398    where

399

400    $$I_M = \frac{N_H\left(\mu_H + \alpha_H + \gamma_H\right)\dfrac{I_H}{N_H}}{ab\left(1 - \left(1 + \dfrac{\gamma_H}{\mu_H}\right)\dfrac{I_H}{N_H}\right)} \tag{36}$$

401

402

403    Before we analyze the sensitivity of the force of infection to the model's parameters

404    related to control, we first deduce a relationship between $\lambda$ and $R_0$.

405    We begin by substituting $I_M$ of equation (36) in equation (35) to obtain

406    $$\lambda = \frac{\left(\mu_H + \alpha_H + \gamma_H\right)\dfrac{I_H}{N_H}}{\left(1 - \left(1 + \dfrac{\gamma_H}{\mu_H}\right)\dfrac{I_H}{N_H}\right)} \tag{37}$$

407

408    If $\alpha_H \approx 0$, the human prevalence, $\frac{I_H}{N_H}$, can be expressed in terms of $R_0$ as follows:

409

410    $$\frac{I_H}{N_H} = \frac{\mu_M\left(1 - g\right)\left(R_0 - 1\right)}{\mu_M\left(1 - g\right)R_0\left(1 + \dfrac{\gamma_H}{\mu_H}\right) + \mu_H ac} \tag{38}$$

411

412

413    Therefore:

414

415

416    $$\lambda = \frac{\mu_M\left(1 - g\right)\left(\mu_H + \alpha_H + \gamma_H\right)\mu_H\left(R_0 - 1\right)}{\mu_M\left(1 - g\right)\left(\mu_H + \gamma_H\right) + \mu_H ac} \tag{39}$$

417
418

419    The partial derivatives of $\lambda$ and $\dfrac{I_H}{N_H}$ with respect to the parameters $\theta_j$ are readily

420    calculated and the sensitivity of $\lambda$ and $\dfrac{I_H}{N_H}$ to the parameters estimated.

421
422
423
424    **Results**
425
426    *Numerical simulations*
427
428         We simulated model (1) with the parameter values available from the literature.

429    However, it is known that these parameters vary with the place, local temperature,

430    climatic factors, mosquito strains and human demography. Therefore, we applied a

431    Monte Carlo simulation algorithm [28] to generate parameter distributions that could

432    mimic real conditions. We used a Beta-distributed random number generator with equal

433    parameters to guarantee the symmetry of the distribution around the mean. Because the

434    Beta distribution with equal parameters has a mean of 0.5, we multiplied the final result

435    by 2. We ran the Monte Carlo algorithm one thousand times to generate the

436    distributions of the parameters. The parameters' baseline values, the mean values of the

437    simulation, the variance and the 95% confidence intervals for each parameter are shown

438    in Table 2.

439
440    *Results of the sensitivity analysis*
441

442         Table 3 shows the results of the sensitivity analysis according to the general

443    equation (27). The results represent the relative amount of variation (expressed in

444    percentual variation)  in the variable if we vary the parameters by 1%.

445

446         Note from Table 3 that $R_0$, $\lambda$ and $\dfrac{I_H}{N_H}$ show the greatest sensitivities to the

447    mosquito's mortality rate $\mu_M$, followed by the biting rate $a$ and the carrying capacity of

448    the immature stages $\kappa_E$. In addition, $R_0$, $\lambda$ and $\dfrac{I_H}{N_H}$ are very insensitive to the larval

449    mortality rate $\mu_E$. Accordingly, a reduction of 1% in the biting rate $a$ or the carrying

450      capacity of the immature stages $\kappa_E$ decreases $R_0$ by 1.94% and 0.69%, respectively, it

451      decreases $\lambda$ by 5.02% and 2.32%, and decreases $\dfrac{I_H}{N_H}$ by 2.67% and 1.34%

452      respectively. Also, an increase of 1% in the mosquito mortality rate $\mu_M$ causes a

453      decrease of 2.42% in $R_0$, of 5.40% in $\lambda$ and of 3.20% in $\dfrac{I_H}{N_H}$. In contrast, increasing

454      the larval mortality rate $\mu_E$ by 1% decreases $R_0$, $\lambda$ and $\dfrac{I_H}{N_H}$ by only 0.000828%,

455      0.00193%, 0.0231% respectively. These differences in the sensitivity of $R_0$, $\lambda$ and

456      $\dfrac{I_H}{N_H}$ to parameter variation can be understood from equation (27). Although the partial

457      derivatives of $\lambda$ with respect to the parameters are smaller than the partial derivatives of

458      $R_0$ with respect to the parameters, the ratio $\dfrac{\theta_j}{\lambda} >> \dfrac{\theta_j}{R_0}$. The same applies for $\dfrac{I_H}{N_H}$.

459

460      **Discussion**

461

462          The knowledge of dengue epidemiology accumulated over the past decades

463      enables us to conclude that the transmission thresholds and the intensity of dengue

464      transmission are determined by several factors: the level of immune protection of the

465      population involved; the serotype of dengue virus circulating at each time; the density,

466      longevity and biting behavior of the mosquitoes; the climate; and the demography of the

467      human hosts [29]. Despite the current development of a safe and effective tetravalent

468      vaccine [1], vector control is still the only available strategy to minimize the number of

469      cases within the affected populations. To date, however, the effectiveness of the

470      strategies for controlling Aedes mosquitoes has been limited. The analysis presented in

471      this paper is intended to contribute to the efforts to check the advance of dengue to areas

472      still free from the disease and to reduce transmission in endemic areas.

473          This paper presents the most complete analysis of what is a basic model for

474      dengue transmission. All the relevant stages are included and the ones not included (like

475      larvae and pupae) can be trivially added to the model.

476          The current paper is an analysis of the basic model proposed in [21] and

477      numerically studied in [20]. The fact that the extrinsic incubation period is changed

478    from being modeled as a fixed time delay to being modeled as an exponentially

479    distributed time period is not relevant for the proposed analysis. As mentioned above,

480    the expressions developed below in this paper with equations (1) can be replaced by the

481    corresponding expressions of the classical Ross-Macdonald model described above by

482    replacing $\dfrac{\gamma_M}{\gamma_M + \mu_M}$ by $e^{-\mu_M \tau}$. In other words, the results of the analysis are the same,

483    irrespective of the way we choose to model the incubation period. Actually, the main

484    difference between this paper and the previous ones [15, 20, 21] is that in the current

485    study we analyze the sensitivity of the endemic equilibrium to variation in the

486    parameters related to transmission in a much more complete way than before. The

487    sensitivity analysis presented in the previous papers consisted only in the derivation of

488    the partial derivatives of $R_0$ with respect to the parameters. This is only part of the

489    sensitivity analysis. In the present paper, the calculation of sensitivity of $R_0$ to the

490    parameters is completed (equation (27)). In addition, we calculated the equilibrium

491    prevalence for the model, obtaining expressions that are completely new, like equation

492    (37) which relates the force of infection to the prevalence of the disease in humans and

493    to the parameters of transmission relative to the humans hosts only. With this

494    expression we propose the estimation of the force of infection for dengue as a function

495    of the equilibrium prevalence in humans.

496    Furthermore, the current and complete sensitivity analysis includes the force of

497    infection and the prevalence of dengue in humans. Finally, the sensitivity of the basic

498    reproduction number and the force of infection to the biting rate is also a quite new

499    finding.

500    Ellis et al. [30] have approached the problem of the sensitivity of dengue by

501    numerically simulating two coupled models, one describing the vector population and

502    the other the host population. These models are extremely complex, including a total of

503    99 parameters for the vector and host populations. Although the calculations based on

504    these models are very important, they mask the dynamics involved. In contrast, the

505    dynamics of dengue constitute the main interest of our paper. Our model contains only

506    16 parameters and admits an analytical solution that can be compared with the classical

507    models designed for the study of vector-borne infections. These differences

508    notwithstanding, the results of Ellis et al. [30] are qualitatively similar to the results that

509    we obtained.

510   Some of the findings of the current paper are qualitatively similar to previous
511   results. However, this is the first paper that proposes a quantification of the relative
512   efficacy of different control strategies. In other words, we are now able to say how
513   much killing adult mosquitoes is more efficient than killing immature stages, for
514   instance.

515   Our results identify the control of adult mosquitoes as the most effective strategy
516   to reduce both $R_0$, $\lambda$ and $\frac{I_H}{N_H}$. However, we are aware that the effectiveness of this
517   strategy is severely constrained, e.g., by the difficulty of achieving sufficiently high
518   coverage of the surfaces used by the mosquitoes for resting [29,31] and by the
519   limitations of ultra-low volume insecticide spraying, which involves a low probability
520   of contact between adult mosquitoes and the insecticide droplets [18].

521   The second most effective strategy is the reduction of the contact between the
522   vectors and hosts, quantified by the daily biting rate $a$. This strategy has been
523   successfully applied in malaria control, e.g., through the use of insecticide-impregnated
524   bed-nets. This approach to malaria control is effective [38] because the malaria
525   mosquito bites at night. Aedes mosquitoes, in contrast, are day-biting mosquitoes, and
526   bed nets are not a feasible method to avoid their bites. In certain countries, however,
527   people habitually take a *siesta*, a rest during the afternoon [18]. In addition, insecticide-
528   treated clothes (ITCs) used as personal protection against malaria infection [18] are
529   beginning to be tested against dengue [10].

530   The next strategy suggested by the analysis of the model's sensitivity involves
531   the carrying capacity of the immature stages, $\kappa_E$. This strategy is associated with the
532   mechanical control of the sources of the mosquitoes. Our assumption is that by
533   destroying mosquitoes' breeding places, we are reducing $\kappa_E$.

534   It is probable that this approach is the most widespread strategy for the control
535   of dengue in endemic regions. However, the results obtained from this strategy have
536   been disappointing. It is probable that these disappointing results are due to the lack of
537   cooperation by the affected communities, which often hampers the application of the
538   method. Unfortunately, $R_0$ was not found to be very sensitive to this strategy. A 1%
539   reduction in $\kappa_E$ yielded only a 0.69% reduction in $R_0$. The force of infection, in
540   contrast, was shown to be relatively sensitive to variation in $\kappa_E$. A 1% reduction in this

541 parameter yielded a 2.32% reduction in $\lambda$. Finally, a 1% reduction in $\kappa_E$ caused a

542 reduction of 1.34% in the human prevalence.

543

544 The least effective strategy analyzed was the use of larvicide. This strategy is

545 expected to increase the mortality rate of immature stages, $\mu_E$. Both $R_0$ and $\lambda$ vary by a

546 fraction on the order of $10^{-3}$ percent, and $\dfrac{I_H}{N_H}$ varies by a fraction on the order of $10^{-2}$

547 percent if we vary $\mu_E$ by 1%.

548 Obviously, the possible control strategies analyzed in this paper are expected to

549 be applied in combination, although we studied each of them in isolation. In addition, it

550 is necessary to carry out a study of financial costs and logistic feasibility to determine

551 the most effective vector control strategy against dengue.

552 The theoretical case of 100% vertical transmission ($g = 1$), i.e., the case in which

553 all of the eggs from the latent and infected mosquitoes are infected, is interesting. In

554 fact, a structural change occurs in our model if $g \rightarrow 1$. The populations of susceptible

555 and infected eggs become completely decoupled. It can be verified that the disease can

556 sustain itself even without human hosts. Actually, as shown by previous authors [32],

557 this is the only way in which the infection circulates exclusively among the vectors in

558 the absence of hosts.

559 In addition, if $g = 1$ and human hosts are introduced into the system, the

560 evolution of the system over time results in a situation in which all mosquitoes are

561 infected because all of the eggs of the infected mosquitoes are infected. Therefore, if

562 $g = 1$ and human hosts are introduced, the population of susceptible mosquitoes and

563 eggs decreases to zero. This result can be verified from equations (8) and (11).

564 Our approach has some important simplifications with respect to reality. The

565 first one is the homogeneously mixing assumption. According to this assumption, the

566 density of every subpopulation is the same everywhere and from the model it seems as

567 if every single infected mosquito has the same probability of contacting every host.

568 Actually, this is not true and it is a notational artifact. In the appendix we explain how

569 this notational artifact can be eliminated. Furthermore, we show how to relax the

570 homogeneously mixing assumption and analyze some consequences of this.

571 The second limitation is that the model predicts a stable endemic equilibrium,

572 which is seldom observed. One reason for this is that in this model, for simplification,

573 we exclude seasonality, which precludes the existence of such equilibrium for long
574 periods of time. The relative sensitivity of the variables to the parameters, however, is
575 also valid (actually to a very good approximation) for non-equilibrium situations. This
576 have already been demonstrated by numerical simulations of a model very similar to the
577 one we are dealing with in this paper [15, 20, 21]. Finally, the actual values of some of
578 the parameters used in the simulations are not known and we had to take advantage of
579 Monte Carlo simulations. The relative sensitivity to the parameters, however, is not
580 affected by the uncertainties in the parameter's values.

581

582

583 **Acknowledgments:** The research from which these results were obtained has received
584 funding from the European Union's Seventh Framework Programme (FP7/2007-2013)
585 under grant agreement no. 282589, from LIM01 HCFMUSP and CNPq. The funders
586 had no role in study design, data collection and analysis, decision to publish, or
587 preparation of the manuscript.

588

589 **Conflicts of Interest:** The authors have declared that no competing interests exist

590

591
592

**Appendix - Some comments on the meaning of the model's equations**

594

595       In this appendix we show how to include spatial heterogeneities in the model

596   and, by doing so, we clarify the meaning of the model's equations.

597       First we assume that mosquitoes have a limited range of flight, which implies

598   that the probability of transmission of infection from one infected mosquito to one

599   susceptible host varies according to the distance between them.

600       Consider the first equation of system (1):

601

602
$$\frac{dS_H}{dt} = -abI_M \frac{S_H}{N_H} - \mu_H S_H + r_H N_H \left(1 - \frac{N_H}{\kappa_H}\right) \tag{A1}$$

603

604   All the variables are <u>densities</u>. This implies that we are considering a very large region

605   where the populations of mosquitoes and hosts are constant, that is, do not vary from

606   point to point. Then, one might think that in equation (A1) a mosquito in a certain place

607   can bite a host which can be very far from it. This is not reasonable and it is not true for

608   equation (A1). To see this, consider the parameter $a$, the mosquitoes' biting rate. We

609   can write this as $a = a'A$, where $a'$ is the biting rate per unit area and $A$ is the area

610   where the mosquitoes' flight ranges. Therefore, only humans inside this area are bitten

611   by this mosquito. But, since the humans and mosquitoes populations are assumed as

612   homogeneously distributed, this does not appear in the equations because in parameter

613   $a$ this effect is hidden.

614    Let us now introduce spatial heterogeneity.  For this we should specify the position $\vec{r}$,

615   representing the spatial location of individuals. Thus, let $S_H(\vec{r})ds$ be the number of

616   human susceptibles in the small area $ds$ around the position $\vec{r}$.

617   Let us now consider how $S_H(\vec{r})ds$ varies with time. Let $I_M(\vec{r}')ds'$ be the number of

618   infected mosquitoes in the small area $ds'$ around the position $\vec{r}'$. The total number of

619   bites the infected mosquitoes population inflicts in a time interval $dt$ is $a'I_M(\vec{r}')ds'dt$. A

620   fraction of those bites $F\left(\left|\vec{r} - \vec{r}'\right|\right)$ is inflicted on the hosts at position $\vec{r}$, that is $S_H(\vec{r})ds$.

621   Of course $F\left(\left|\vec{r} - \vec{r}'\right|\right)$ is a decreasing function of the distance $\left|\vec{r} - \vec{r}'\right|$ between infected

622   mosquitoes and susceptible humans. Thus equation (A1) becomes:

623

$$\frac{dS_H(\vec{r})}{dt} = -b\frac{S_H(\vec{r})}{N_H(\vec{r})}\int d\vec{s}' a'(\vec{r}') F(|\vec{r}-\vec{r}'|) I_M(\vec{r}') - \mu_H S_H(\vec{r}) + r_H N_H(\vec{r})\left(1 - \frac{N_H(\vec{r})}{\kappa_H(\vec{r})}\right) \text{(A2)}$$

625

626  All the other equations in system (1) should be similarly modified and, of course, the

627  result is very difficult to integrate. When $a'(\vec{r}') F(|\vec{r}-\vec{r}'|)$ is equal to $a'A\theta(|\vec{r}-\vec{r}'|)$, and

628  the densities are homogenously distributed in space, we have

629

630  $$b\frac{S_H(\vec{r})}{N_H(\vec{r})}\int d\vec{s}' a'(\vec{r}') F(|\vec{r}-\vec{r}'|) = b\frac{S_H}{N_H} I_M \int d\vec{s}' a'(\vec{r}') F(|\vec{r}-\vec{r}'|) = b\frac{S_H}{N_H} I_M a \qquad \text{(A3)}$$

631

632  and equation (A2) reduces to (A1).

633

634  The above formalism is necessary when we are dealing with large regions of space,

635  where heterogeneities are significant. However, for small regions, where heterogeneities

636  can be neglected, the system of equations (1) of the main text, are a good

637  approximation. The relative sensitivity of the transmission variables to the studied

638  parameters, however, is not expected to be significantly influenced by spatial

639  heterogeneities. Of course, the value of the transmission variables may vary from place

640  to place but the <u>relative</u> sensitivity, the main objective of the present analysis, of these

641  variables to the parameters should be the same.

642 **References**

643

644 1. Guy B, Almond J, Lang J (2011) Dengue vaccine prospects. Lancet 377: 381–382.

645 2. Gubler DJ (2002) The global emergence/resurgence of arboviral diseases as public
646 health problems. Arch Med Res 33: 330–342.

647 3. WHO (2012) Dengue and severe dengue. Fact sheet N°117. Available:
648 http://www.who.int/mediacentre/factsheets/fs117/en/. Accessed 10 April 2012.

649 4. WHO (2009) Dengue and dengue haemorrhagic fever. Fact sheet Nº117. Available:
650 http://who.int/mediacentre/factsheets/fs117/en/print.html. Accessed 11 March 2011.

651 5. Gubler DJ (2002) The global emergence/resurgence of arboviral diseases as public
652 health problems. Arch Med Res 33: 330–342.

653 6. Beatty ME, Letson GW, Margolis HS (2008) Estimating the global burden of dengue.
654 Abstract book: dengue 2008. Proceedings of the 2nd International Conference on
655 Dengue and Dengue Haemorrhagic Fever, Phuket, Thailand.

656 7. World Tourism Organization (2011) Tourism highlights. Available: www.world-
657 tourism.org/facts/menu.html. Accessed 11 March 2011.

658 8. Suaya JA, Shepard DS, Siqueira JB (2009) Cost of dengue cases in eight coutries in
659 the Americas and Asia: a propective study. Am J Trop Med Hyg 80: 846–855.

660 9. Khasnis AA, Nettleman MD (2005) Global warming and infectious disease. Arch
661 Med Res 36: 689–696.

662 10. Wilder-Smith A, Renhorn KE, Tissera H, Abu Bakar S, Alphey L, et al. (2012)
663 DengueTools: innovative tools and strategies for the surveillance and control of dengue.
664 Glob Health Action 5: 17273.

665 11. Gubler DJ (2011) Dengue, urbanization and globalization: the unholy trinity of the
666 21$^{st}$ century. Trop Med Health 39(4 Suppl): 3–11.

667 12. Beatty ME, Beutels P, Meltzer MI, Shepard DS, Hombach J, et al. (2011) Health
668 economics of dengue: a systematic literature review and expert panel's assessment. Am
669 J Trop Med Hyg 84: 473–488.

670    13. Wilder-Smith A, Ooi EE, Vasudevan SG, Gubler DJ (2010) Update on dengue:
671    epidemiology, virus evolution, antiviral drugs, and vaccine development. Curr Infect
672    Dis Rep 12: 157–164.

673    14. Lambrechts L, Scott TW, Gubler DJ (2010) Consequences of the expanding global
674    distribution of Aedes albopictus for dengue virus transmission. PLoS Negl Trop Dis 4:
675    e646.

676    15. Massad E, Coutinho FAB, Lopez LF, da Silva DR (2011) Modeling the impact of
677    global warming on vector-borne infections. Phys Life Rev 8: 169–199.

678    16. Luz PM, Vanni T, Medlock J, Paltiel AD, Galvani AP (2011) Dengue vector control
679    strategies in an urban setting: an economic modelling assessment. Lancet 377: 1673–
680    1680.

681    17. Massad E, Coutinho FAB (2011) The cost of dengue control. Lancet 377: 1630–
682    1631.

683    18. Reiter P, Gubler DJ (2001) Surveillance and control of urban dengue vectors. In:
684    Gubler DJ, Kuno G, editors. Dengue and dengue hemorrhagic fever. Wallingford, UK:
685    CABI Publishing. pp. 425–462.

686    19. MacDonald G (1952) The analysis of equilibrium in malaria. Trop Dis Bull 49:
687    813–828.

688    20. Burattini MN, Chen M, Chow A, Coutinho FAB, Goh KT, et al. (2008) Modelling
689    the control strategies against dengue in Singapore. Epidemiol Infect 136: 309–319.
690

691    21. Coutinho FAB, Burattini MN, Lopez LF, Massad E (2006) Threshold conditions for
692    a non-autonomous epidemic system describing the population dynamics of dengue. Bull
693    Math Biol 68: 2263–2282.
694

695    22. Coutinho FAB, Burattini MN, Lopez LF, Massad E (2005) An approximate
696    threshold condition for non-autonomous system: an application to a vector-borne
697    infection. Math Comp Simul 70: 149–158.

698   23. Massad E, Behrens RH, Burattini MN, Coutinho FAB (2009) Modeling the risk of
699       malaria for travelers to areas with stable malaria transmission. Malar J 8: 296.

700   24. Chitnis N, Hyman JM, Cushing JM (2008) Determining important parameters in the
701       spread of malaria through the sensitivity analysis of a mathematical model. Bull Math
702       Biol 70: 1272–1296.

703   25. Lopez LF, Coutinho FAB, Burattini MN, Massad E (2002) Threshold conditions for
704       infection persistence in complex host-vectors interactions. C R Biol 325: 1073–1084.

705   26. Fegan G, Noor AM, Akhwale WS, Cousens S, Snow RW (2007) Effect of expanded
706       insecticide-treated bednet coverage on child survival in rural Kenya: a longitudinal
707       study. Lancet 370: 1035–1039.

708   27. Ross R (1911) The prevention of malaria, 2$^{nd}$ ed., with addendum on the theory of
709       happenings. London: Murray.

710   28. Amaku M, Azevedo RS, Castro RM, Massad E, Coutinho FAB (2009) Relationship
711       among epidemiological parameters in a non-immunized Brazilian community. Mem
712       Inst Oswaldo Cruz 104: 897–900.

713   29. Rodrigues HS, Monteiro MT, Torres DFM (2012). Dengue in Cape Verde: vector
714       control and vaccination. arXiv: 1204.0544v1.

715   30. Ellis AM, Garcia AJ, Focks DA, Morrison AC, Scott TW (2011) Parameterization
716       and sensitivity analysis of a complex simulation model for mosquito population
717       dynamics, dengue transmission and their control. Am J Trop Med Hyg 82: 257–264.

718   31. Integrated Vector Management (2012) Available:
719       http://www.ivmproject.net/about/index.cfm?fuseaction=static&label=dengue. Accessed
720       1 April 2012.
721

722   32. Adams B, Boots M (2010) How important is vertical transmission in mosquitoes for
723       the persistence of dengue? Insights from a mathematical model. Epidemics 2: 1–10.

724   33. Yasuno M, Tonn RJ (1970) A study of biting habits of Aedes aegypti in Bangkok,
725       Thailand. Bull World Health Organ 43: 319–325.

726    34. Ocampo CB, Wesson DM (2004) Population dynamics of Aedes aegypti from a

727    dengue hyperendemic urban setting in Colombia. Am J Trop Med Hyg **7**1: 506–513.

728    35. Index Mundi (2011) Available: http://www.indexmundi.com/map/?v=30&l=pt.

729    Accessed 18 August 2011.

730

731    35. Halstead SB (1990) Dengue. In: Warren KS, Mahmoud AAF editors. Tropical and

732    geographical medicine. New York: McGraw-Hill. pp. 675–684.

733    37. Forattini OP (1996) Medical culicidology. São Paulo: EDUSP.

734    38. Brownstein JS, Heth E, O'Neill L (2003) The potential of virulent Wolbachia to

735    modulate disease transmission by insects. J Invertebr Pathol 84**:** 24–29.

736

737

**Table 1.** Model variables and their biological meanings.

739

| Variable | Biological Meaning |
|---|---|
| $S_H$ | Density of susceptible humans |
| $I_H$ | Density of infected humans |
| $R_H$ | Density of recovered humans |
| $S_M$ | Density of uninfected mosquitoes |
| $L_M$ | Density of latent mosquitoes |
| $I_M$ | Density of infected mosquitoes |
| $S_E$ | Density of uninfected eggs (imm. Stages) |
| $I_E$ | Density of infected aquatic forms |

740

741

**Table 2.** Model parameters, biological meaning, values and sources. The mean, variance and 95% CI were obtained with Monte Carlo simulations.

| Parameter | Meaning | Value (Baseline) | Mean | Variance | 95% CI | Source |
|---|---|---|---|---|---|---|
| $a$ | Average daily rate of biting | 0.164 | 0.1682 | 0.026 | $9.8 \times 10^{-3}$ | [33] |
| $b$ | Fraction of bites actually infective | 0.6 | 0.6062 | 0.296 | 0.0337 | [34] |
| $\mu_H$ | Human natural mortality rate | $3.5 \times 10^{-5}$ days$^{-1}$ | $3.55 \times 10^{-5}$ | $1.019 \times 10^{-9}$ | $2.00 \times 10^{-6}$ | [35] |
| $r_H$ | Birth rate of humans | $9.5 \times 10^{-5}$ days$^{-1}$ | $9.531 \times 10^{-5}$ | $8.959 \times 10^{-9}$ | $5.3 \times 10^{-6}$ | [35] |
| $\kappa_H$ | Carrying capacity of humans | $5 \times 10^{6}$ | $5.0123 \times 10^{6}$ | $2.052 \times 10^{13}$ | $2.81 \times 10^{5}$ | [35] |
| $\alpha_H$ | Dengue mortality in humans | $3.5 \times 10^{-4}$ days$^{-1}$ | $3.473 \times 10^{-4}$ | $1.00 \times 10^{-7}$ | $1.97 \times 10^{-5}$ | [36] |
| $\gamma_H$ | Human recovery rate | 0.143 days$^{-1}$ | 0.1434 | 0.017 | $8.097 \times 10^{-3}$ | [36] |
| $P$ | Hatching rate of susceptible eggs | 0.15 days$^{-1}$ | 0.151 | 0.019 | $8.55 \times 10^{-3}$ | [37] |
| $\gamma_M$ | Latency rate in mosquitoes | 0.143 days$^{-1}$ | 0.1434 | 0.017 | $8.097 \times 10^{-3}$ | [20] |
| $\mu_M$ | Natural mortality rate of mosquitoes | 0.09 days$^{-1}$ | 0.08329 | $1.5 \times 10^{-4}$ | $5.52 \times 10^{-3}$ | [38] |
| $r_M$ | Oviposition rate | 50 days$^{-1}$ | 51.8295 | 2073.9 | 2.8226 | [38] |
| $G$ | Proportion of infected eggs | 0.1 | 0.0964 | 0.008 | $5.684 \times 10^{-3}$ | Assumed |
| $\kappa_E$ | Carrying capacity of eggs | $9.8 \times 10^{7}$ | $9.787 \times 10^{7}$ | $8.003 \times 10^{15}$ | $5.545 \times 10^{6}$ | Assumed |
| $\mu_E$ | Natural mortality rate of eggs | 0.1 days$^{-1}$ | 0.101 | 0.008 | $5.6644 \times 10^{-3}$ | [38] |
| $C$ | Dengue susceptibility of A. aegypti | 0.54 | 0.5265 | 0.249 | 0.03191 | [34] |
| $c_S$ | Climatic factor | 0.07 | 0.07 | 0.004 | 0.00398 | Assumed |

**Table 3.** Results of the sensitivity analysis according to the general equation (27). The results represent the relative amount of variation (expressed in percentual variation) in the variable if we vary the parameters by 1%.

| Variable | Mean | 95% Confidence Interval |
|---|---|---|
| $R_0$ | 1.74 | $1.45 - 2.07$ |
| $\lambda$ | $2.59 \times 10^{-5}$ | $1.48 \times 10^{-5} - 3.96 \times 10^{-5}$ |
| $I_H / N_H$ | $1.04 \times 10^{-4}$ | $3.84 \times 10^{-5} - 1.34 \times 10^{-4}$ |

### Sensitivity of $R_0$ to the control parameters

| Parameter | Mean |
|---|---|
| $a$ | 1.94 |
| $\kappa_E$ | 0.69 |
| $\mu_E$ | (-) $8.28 \times 10^{-4}$ |
| $\mu_M$ | (-) 2.42 |

### Sensitivity of $\lambda$ to the control parameters

| Parameter | Mean |
|---|---|
| $a$ | 5.02 |
| $\kappa_E$ | 2.32 |
| $\mu_E$ | (-) $1.93 \times 10^{-3}$ |
| $\mu_M$ | (-) 5.40 |

### Sensitivity of $I_H / N_H$ to the control parameters

| Parameter | Mean |
|---|---|
| $a$ | 2.67 |
| $\kappa_E$ | 1.34 |
| $\mu_E$ | (-) $2.31 \times 10^{-2}$ |
| $\mu_M$ | (-) 3.20 |